\documentclass[preprint,10pt]{elsarticle}

\usepackage{mathrsfs,amsmath,amssymb}

\usepackage{mathrsfs,amsmath,amssymb}
\usepackage{mathrsfs,amsmath,amssymb}
\usepackage{amscd} 
\usepackage{relsize,amsmath} 

\makeatletter
\newcommand{\biggg}[1]{{\hbox{$\left#1\vbox to 20.5pt{}\right.\n@space$}}}
\newcommand{\Biggg}[1]{{\hbox{$\left#1\vbox to 23.5pt{}\right.\n@space$}}}
\newcommand{\bigggg}[1]{{\hbox{$\left#1\vbox to 26.5pt{}\right.\n@space$}}}
\newcommand{\Bigggg}[1]{{\hbox{$\left#1\vbox to 29.5pt{}\right.\n@space$}}}
\newcommand{\biggggg}[1]{{\hbox{$\left#1\vbox to 32.5pt{}\right.\n@space$}}}
\newcommand{\Biggggg}[1]{{\hbox{$\left#1\vbox to 35.5pt{}\right.\n@space$}}}
\newcommand{\bigggggg}[1]{{\hbox{$\left#1\vbox to 38.5pt{}\right.\n@space$}}}
\newcommand{\Bigggggg}[1]{{\hbox{$\left#1\vbox to 41.5pt{}\right.\n@space$}}}
\makeatother

\makeatletter
\@addtoreset{equation}{section}
\makeatother

\usepackage{braket}

\begin{document}

\newtheorem{thm}{Theorem}
\newtheorem{lem}[thm]{Lemma}
\newdefinition{rmk}{Remark}
\newproof{pf}{Proof}
\newproof{pot}{Proof of Theorem \ref{thm2}}

\begin{frontmatter}



\title{Asymptotic behavior of solutions toward 
the rarefaction waves 
to the Cauchy problem for 
the generalized Benjamin-Bona-Mahony-Burgers equation with dissipative term}


\author[labe1
]{Natsumi Yoshida}
\ead{14v00067@gmail.com
}

\address[label]{Graduate Faculty of Interdisciplinary Research Faculty of Education, 
University of Yamanashi, Kofu, Yamanashi 400-8510, Japan.}
\address{}

\begin{abstract}
In this paper, we investigate the asymptotic behavior of solutions 
to the Cauchy problem with the far field condisiton 
for the generalized Benjamin-Bona-Mahony-Burgers equation 
with a fourth-order dissipative term. 
When the corresponding Riemann problem 
for the hyperbolic part 
admits a Riemann solution which 
consists of single rarefaction wave, 
it is proved that 
the solution of the Cauchy problem tends toward the rarefaction wave 
as time goes to infinity. 
We can further obtain the same global asymptotic stability 
of the rarefaction wave 
to the generalized Korteweg-de Vries-Benjamin-Bona-Mahony-Burgers equation 
with a fourth-order dissipative term 
as the former one. 
\end{abstract}

\begin{keyword} 
Benjamin-Bona-Mahony-Burgers equation 
\sep Korteweg-de Vries-Benjamin-Bona-Mahony-Burgers equation \sep convex flux 
\sep asymptotic behavior \sep rarefaction wave 

\medskip
AMS subject classifications: 35K55, 35B40, 35L65
\end{keyword}

\end{frontmatter}

%



\pagestyle{myheadings}
\thispagestyle{plain}
\markboth{N. YOSHIDA}{GENERALIZED BENJAMIN-BONA-MAHONY-BURGERS EQUATION}

\section{Introduction and main theorems}
In this paper, 
we consider the asymptotic behavior of solutions to the Cauchy problem 
for the generalized Benjamin-Bona-Mahony-Burgers equation 
with a fourth-order dissipative term 
\begin{eqnarray}
 \left\{\begin{array}{ll}
  \partial_tu +\partial_x \big( \, 
f(u) \, \big) - \alpha \, \partial_t\partial_x^2 u 
- \beta \, \partial_x^2 u 
+\gamma  \, \partial_x^4 u=0
  \qquad &\big( \, t>0, \: x\in \mathbb{R} \, \big), \\[5pt]
 u(0, x) = u_0(x) \rightarrow u_{\pm} \qquad &( x \rightarrow \pm \infty ).
 \end{array}
 \right.\,
\end{eqnarray}
Here, $u=u(t, x)$ is 
the unknown function of $t>0$ and $x\in \mathbb{R}$, 
and $\alpha$, $\beta$, $\gamma$ are positive constants
$u_0$ is the initial data, 
and $u_{\pm } \in \mathbb{R}$ 
are the prescribed far field states. 
We suppose that $f$ is a smooth function. 

There are many results concerning with 
the mathematical structure, 
such as the global existence and time-decay properties of solutions, 
of the generalized Benjamin-Bona-Mahony-Burgers equation 
with the dissipative terms (see 
Kondo-Webler \cite{kondo-webler1}, \cite{kondo-webler2}, 
\cite{kondo-webler3}, \cite{kondo-webler4}, 
Wang \cite{wang}, Xu-Li \cite{xu-li}
Zhao-Xuan \cite{zhao-xuan} 
and so on). 
The model (1.1) is closely related to the following Benjamin-Bona-Mahony-Burgers equation 
\begin{equation}
\partial_tu + \delta \, \partial_x u
+u \, \partial_x u - \alpha \, \partial_t\partial_x^2 u 
- \beta \, \partial_x^2 u=0,
\end{equation}
where $\delta \in \mathbb{R}$. 
The mathematical structure of (1.2) 
have also been investigated by Amick-Bona-Schonbek \cite{amik-bona-schonbek}, 
Mei \cite{mei1}, \cite{mei2}, Mei-Schmeiser \cite{mei-schmeiser},
Naumkin \cite{naumkin} and so on. 
When $\beta=0$, then (1.2) becomes the following Benjamin-Bona-Mahony equation 
\begin{equation}
\partial_tu + \delta \, \partial_x u
+u \, \partial_x u - \alpha \, \partial_t\partial_x^2 u =0,
\end{equation}
which was advocated by \cite{benjamin-bona-mahony} 
as a refinement of the following Korteweg de-Vries equation 
\begin{equation}
\partial_tu + \delta \, \partial_x u
+u \, \partial_x u - \alpha \, \partial_x^3 u =0.
\end{equation}
For the case $\alpha=1$ and $\beta=0$ of (1.2), (1.2) is 
the so-called regularized long wave equation, 
which was proposed 
by Peregrine \cite{peregrine} and \cite{benjamin-bona-mahony}, 
as follows. 
\begin{equation}
\partial_tu + \delta \, \partial_x u
+u \, \partial_x u - \partial_t\partial_x^2 u =0.
\end{equation}
We note that (1.3) and (1.5) are known as the 
approximated models for the long waves of small amplitude.

We are going to obtain the rarefaction stability of the solution to (1.1). 
Therefore we deal with the case 
where the flux function $f$ is fully convex, 
that is, 
\begin{equation}
f''(u)>0\quad (u \in \mathbb{R}), 
\end{equation}
and $u_-<u_+$.
Then, since the corresponding Riemann problem (cf. \cite{lax}) 
\begin{eqnarray}
 \left\{\begin{array} {ll}
 \partial _t u + \partial _x \bigl( \, f(u) \, \bigr)=0 
 \, \: \; \qquad \big( \, t>0,\: x\in \mathbb{R} \, \big),\\[5pt]
u(0, x)=u_0 ^{\rm{R}} (x)
      := \left\{\begin{array} {ll}
         u_-  & \; (x < 0),\\[5pt]
         u_+  & \; (x > 0)
         \end{array}\right.
 \end{array}
  \right.\,
\end{eqnarray}
turns out to admit a single
rarefaction wave solution, 
we expect that 
the solution of the Cauchy problem (1.1) tends toward
the rarefaction wave 
as time goes to infinity. 
Here, the rarefaction wave connecting $u_-$ to $u_+$ 
is given by 
\begin{equation}
u^r \left( \, \frac{x}{t}\: ;\:  u_- ,\:  u_+ \right)
= \left\{
\begin{array}{ll}
  u_-  & \; \bigl(\, x \leq f'(u_-)\,t \, \bigr),\\[7pt]
  \displaystyle{ (f')^{-1}\left( \frac{x}{t}\right) } 
  & \; \bigl(\, f'(u_-)\,t \leq x \leq f'(u_+)\,t\,  \bigr),\\[7pt]
   u_+ & \; \bigl(\, x \geq f'(u_+)\,t \, \bigr).
\end{array}
\right. 
\end{equation} 
In particular, we also expect that if $u_- = u_+=:\tilde{u}$, 
then the solution of the Cauchy problem (1.1) 
tends toward the constant state $\tilde{u}$ 
as time goes to infinity. 

There are many results conserning with the rarefaction stabilities. 
For viscous conservation law, 
\begin{eqnarray*}
 \left\{\begin{array}{ll}
  \partial_tu +\partial_x \bigl( \, 
f(u) - \mu \, \partial_xu 
 \, \bigr) =0
  \qquad &\big( \, t>0, \: x\in \mathbb{R} \, \big), \\[5pt]
 u(0, x) = u_0(x) \rightarrow u_{\pm} \qquad &( x \rightarrow \pm \infty ),
 \end{array}
 \right.\,
\end{eqnarray*}
with the condition (1.6), 
Il'in-Ole{\u\i}nik \cite{ilin-oleinik} showed 
that the solution tends toward the single rarefaction wave under the condition $u_-<u_+$ 
(and the one does the single viscous shock wave under the condition $u_{-}>u_{+}$, 
for further studies, 
see \cite{matsu-nishi3}, \cite{matsumura-yoshida}, \cite{yoshida1}, \cite{yoshida8} and so on). 
Hattori-Nishihara \cite{hattori-nishihara} also obtained the pointwise 
and time-decay estimates of the difference 
$|\, u-u^r \,|$. 
Harabetian \cite{harabetian} further considered the following 
rarefaction problem for a quasilinear parabolic equation 
\begin{eqnarray*}
 \left\{\begin{array}{ll}
  \partial_tu +\partial_x \bigl( \, 
f(u) -A'(u) \, \partial_xu  
 \, \bigr) =0
  \qquad &\big( \, t>0, \: x\in \mathbb{R} \, \big), \\[5pt]
  u(0, x) = u_0(x) \rightarrow u_{\pm} \qquad &( x \rightarrow \pm \infty ),
 \end{array}
 \right.\,
\end{eqnarray*}
where $A'(u)\ge0 $ $(u \in \mathbb{R})$, and obtained the precise time-decay estimates of global stability 
of the rarefaction wave with the aid of 
the arguments on monotone semigroups by Osher-Ralston \cite{osh-ral}. 
For the following Cauchy problem of the Matsumura-Nishihara model
\begin{eqnarray*}
 \left\{\begin{array}{ll}
  \partial_tu +\partial_x \bigl( \, 
f(u) - \mu \, | \,  \partial_xu \, |^{p-1} \, \partial_xu  
 \, \bigr) =0
  \qquad &\big( \, t>0, \: x\in \mathbb{R} \, \big), \\[5pt]
  u(0, x) = u_0(x) \rightarrow u_{\pm} \qquad &( x \rightarrow \pm \infty ),
 \end{array}
 \right.\,
\end{eqnarray*}
where $p>1$ and the viscosity $\mu \, | \,  \partial_xu \, |^{p-1} \, \partial_xu  $ 
is the so-called Ostwald-de Waele-type viscosity 
advocated by de Waele \cite{de waele} and Ostwald \cite{ost}
(which is a typical example for the non-Newtonian viscosity, 
see also \cite{chh1}, \cite{chh2}, \cite{chh-ric}, 
\cite{jah-str-mul}, \cite{lad1}, \cite{ma}, \cite{ma-pr-st}, 
\cite{soc} and so on), 
Matsumura-Nishihara \cite{matsu-nishi2} first investigated and proved the global stability 
of the rarefaction wave 
by using the technical energy method. 
Yoshida \cite{yoshida2} further obtained its precise time-decay estimates 
by using the time-weighted energy method 
(for the stabilities of the multiwave pattern, 
see \cite{yoshida3'}, \cite{yoshida3}, \cite{yoshida4}). 
Furthermore, Matsumura-Yoshida \cite{matsumura-yoshida'} considered 
the following Cauchy problem of the non-Newtonian viscous conservation law  
\begin{eqnarray*}
 \left\{\begin{array}{ll}
  \partial_tu +\partial_x \Big( \, 
f(u) - \sigma \big(\, \partial_xu \, \big)
 \, \Big) =0
  \qquad &\big( \, t>0, \: x\in \mathbb{R} \, \big), \\[5pt]
 u(0, x) = u_0(x) \rightarrow u_{\pm} \qquad &( x \rightarrow \pm \infty ),
 \end{array}
 \right.\,
\end{eqnarray*}
where viscosity function $\sigma$ satisfies the conditions
\begin{eqnarray*}
 \left\{\begin{array}{ll}
\sigma(0)=0, \quad 
\sigma'(v)>0 \quad (v \in \mathbb{R}), \\[5pt]
| \,  \sigma(v) \, |\sim |\, v \, |^p, \quad 
| \,  \sigma'(v) \, |\sim |\, v \, |^{p-1} \quad 
\big( \,  |\, v \, |\rightarrow \infty \, \big),
\end{array} 
\right.\,
\end{eqnarray*}
and obtained the global stability of the rarefaction wave for the case 
$0<p<3/7$. 
Recently, Yoshida \cite{yoshida7} obtained this rarefaction stability for more general case 
$0<p<1/3$ and its precise time-decay estimates. 
For the rarefaction problem of the Korteweg-de Vries equation 
\begin{eqnarray*}
 \left\{\begin{array}{ll}
  \partial_tu +\partial_x \bigl( \, 
-3 \, u^2 + \partial_x^2u 
 \, \bigr) =0
  \qquad \big( \, t>0, \: x\in \mathbb{R} \, \big), \\[5pt]
  u(0, x) = u_0(x) 
  \left\{\begin{array}{ll}
  \rightarrow 0 \quad &(x\rightarrow \infty),\\[5pt]
  \rightarrow u_{-}>0 \quad &(x\rightarrow -\infty),
   \end{array}
 \right.\,
  \end{array}
 \right.\,
\end{eqnarray*}
Egorova-Grunert-Teschl \cite{ego-gru-tes} and Egorova-Teschl \cite{ego-tes} 
obtained the existence and uniqueness in some class of the classical solution 
, and Andreiev-Egorova-Lange-Teschl \cite{and-ego-lan-tes} 
also obtained the valid asymptotic formula of the solution. 
For diffusive dispersive conservation laws, 
Wang-Zhu \cite{wan-zhu} 
obtained 
the local stability of the rarefaction wave for the 
following Cauchy problem 
of the generalized Korteweg-de Vries-Burgers equation 
\begin{eqnarray*}
 \left\{\begin{array}{ll}
  \partial_tu +\partial_x \bigl( \, 
f(u) - \mu \, \partial_xu + \delta \, \partial_x^2u 
 \, \bigr) =0
  \qquad &\big( \, t>0, \: x\in \mathbb{R} \, \big), \\[5pt]
 u(0, x) = u_0(x) \rightarrow u_{\pm} \qquad &( x \rightarrow \pm \infty ).
 \end{array}
 \right.\,
\end{eqnarray*}
Duan-Zhao \cite{dua-zha} and Yoshida \cite{yoshida6} 
further obtained the global stabilities of the rarefaction wave 
under some growth conditions for $f$ 
(for stability and time-decay properties of a travelling wave, see 
\cite{bona-schonbek}, \cite{bona-rajopadhye-schonbek}, \cite{nishi-raj}). 
For the Cauchy problem 
of the generalized Korteweg-de Vries-Burgers-Kuramoto equation 
\begin{eqnarray*}
 \left\{\begin{array}{ll}
  \partial_tu +\partial_x \bigl( \, 
f(u) - \mu \, \partial_xu + \delta \, \partial_x^2u 
+ \nu \, \partial_x^3u 
 \, \bigr) =0
  \qquad &\big( \, t>0, \: x\in \mathbb{R} \, \big), \\[5pt]
 u(0, x) = u_0(x) \rightarrow u_{\pm} \qquad &( x \rightarrow \pm \infty ), 
 \end{array}
 \right.\,
\end{eqnarray*}
Ruan-Gao-Chen \cite{rua-gao-che} first obtained the local stability 
of the rarefaction wave. 
Duan-Fan-Kim-Xie \cite{dua-fan-kim-xie} and Yoshida \cite{yoshida6} 
also obtained the global stability 
of the rarefaction wave 
under some growth conditions for $f$.
Recently, Yoshida \cite{yoshida10} (see also \cite{yoshida9}) 
investigated the following Cauchy problem for 
a diffusive dispersive conservation waw, that is, 
the generalized Korteweg-de Vries-Burgers-Kuramoto equation without the viscosity term as 
\begin{eqnarray*}
 \left\{\begin{array}{ll}
  \partial_tu +\partial_x \bigl( \, 
f(u) + \delta \, \partial_x^2u 
+ \nu \, \partial_x^3u 
 \, \bigr) =0
  \qquad &\big( \, t>0, \: x\in \mathbb{R} \, \big), \\[5pt]
 u(0, x) = u_0(x) \rightarrow u_{\pm} \qquad &( x \rightarrow \pm \infty ), 
 \end{array}
 \right.\,
\end{eqnarray*}
and obtained the global stability of the rarefaction wave 
under some growth conditions for $f$. 


Our main results of the present paper are as follows. 

\medskip

\noindent
{\bf Theorem 1.1} (Main Theorem I){\bf .}\quad{\it
Assume the far field states $u_{\pm}$ satisfy $u_- = u_+=\tilde{u}$, and 
the convective flux $f \in C^2(\mathbb{R})$. 
Further assume the initial data satisfy
$u_0-\tilde{u} \in L^2$ and
$\partial _xu_0 \in H^2$. 
Then the Cauchy problem {\rm(1.1)} has a 
unique global in time 
solution $u$ 
satisfying 
\begin{eqnarray*}
 \left\{\begin{array}{ll}
u-\tilde{u} \in C^0\bigl( \, [\, 0, \, \infty \, ) \, ; H^3 \, \bigr), \\[5pt]
\partial _x u \in L^2\bigl( \, 0,\, \infty \, ; H^3 \, \bigr), \\[5pt]
\partial _t u \in L^2\bigl( \, 0,\, \infty \, ; H^2 \, \bigr),
\end{array} 
\right.\,
\end{eqnarray*}
and the asymptotic behavior 
$$
\lim _{t \to \infty}
\bigg( \,
\sup_{x\in \mathbb{R}}
| \, u(t,x) - \tilde{u} \, | 
+\sup_{x\in \mathbb{R}}
 |\,\partial_xu (t, x)\,| 
+\sup_{x\in \mathbb{R}}
 |\,\partial_x^2u (t, x)\,|
\, \bigg)= 0.
$$
}

\medskip

\noindent
{\bf Theorem 1.2} (Main Theorem I\hspace{-.1em}I){\bf .}\quad{\it
Assume the far field states $u_{\pm}$ satisfy $u_- < u_+$, and 
the convective flux $f \in C^5(\mathbb{R})$ 
satisfy {\rm(1.6)}. 
Further assume the initial data satisfy
$u_0-u_0 ^{\rm{R}} \in L^2$ and
$\partial _xu_0 \in H^2$. 
Then the Cauchy problem {\rm(1.1)} has a 
unique global in time 
solution $u$ 
satisfying 
\begin{eqnarray*}
 \left\{\begin{array}{ll}
u-u_0 ^{\rm{R}} \in C^0\bigl( \, [\, 0, \, \infty \, ) \, ; H^3 \, \bigr), \\[5pt]
\partial _x u \in L^2_{\rm{loc}}\bigl( \, 0,\, \infty \, ; H^3 \, \bigr), \\[5pt]
\partial _t u \in L^2_{\rm{loc}}\bigl( \, 0,\, \infty \, ; H^2 \, \bigr),
\end{array} 
\right.\,
\end{eqnarray*}
and the asymptotic behavior 
\begin{eqnarray*}
\left\{\begin{array}{ll}
\displaystyle{
\lim _{t \to \infty}\sup_{x\in \mathbb{R}} \, 
\left| \, 
u(t,x) - u^r\left(\, \frac{x}{t}\: ;\: u_-, \, u_+ \, \right) 
\, \right|} = 0,\\[10pt]
\displaystyle{
\lim _{t \to \infty}\sup_{x\in \mathbb{R}} \, 
| \, 
\partial_{x}u(t,x) 
- \partial_{x}u^r\left(\, t, \, x\: ;\: u_-, \, u_+ \, \right)
\, |} = 0, \\[10pt]
\displaystyle{
\lim _{t \to \infty}\sup_{x\in \mathbb{R}} \, 
| \, 
\partial_{x}^2u(t,x) 
- \partial_{x}^2u^r\left(\, t, \, x\: ;\: u_-, \, u_+ \, \right)
\, |} = 0, 
\end{array} 
\right.\,
\end{eqnarray*}
where $\partial_{x}u^r$ and $\partial_{x}^2u^r$ 
are given by 
\begin{equation*}
\partial_{x}u^r\left(\, t, \, x \: ;\: u_-, \, u_+ \, \right)
= \left\{
\begin{array}{ll}
  0  & \; \bigl(\, x \leq f'(u_-)\,t \, \bigr),\\[7pt]
  \displaystyle{ 
  \frac{1}
  {f''\bigg( \, (f')^{-1}\left( \displaystyle{\frac{x}{t}} \right) \, \bigg)} 
  \, \frac{1}{t}
  } 
  & \; \bigl(\, f'(u_-)\,t \leq x \leq f'(u_+)\,t\,  \bigr),\\[15pt]
  0 & \; \bigl(\, x \geq f'(u_+)\,t \, \bigr),
\end{array}
\right. 
\end{equation*} 
\begin{align*}
\begin{aligned}
&\partial_{x}^2u^r\left(\, t, \, x \: ;\: u_-, \, u_+ \, \right)\\
&= \left\{
\begin{array}{ll}
  0  & \; \bigl(\, x \leq f'(u_-)\,t \, \bigr),\\[7pt]
  \displaystyle{ 
  \frac{1}
  {\left( \, f''\bigg( \, (f')^{-1}\left( \displaystyle{\frac{x}{t}} \right) \, \bigg) \, \right)^3} 
  \, \frac{x}{t^3}
  -\frac{1}
  {f''\bigg( \, (f')^{-1}\left( \displaystyle{\frac{x}{t}} \right) \, \bigg)} 
  \, \frac{1}{t^2}
  } 
  & \; \bigl(\, f'(u_-)\,t \leq x \leq f'(u_+)\,t\,  \bigr),\\[15pt]
  0 & \; \bigl(\, x \geq f'(u_+)\,t \, \bigr).
\end{array}
\right. 
\end{aligned}
\end{align*}
}

\medskip

Furthermore, for the Cauchy problem 
for the generalized Benjamin-Bona-Mahony-Burgers equation 
with third-order dispersive and fourth-order dissipative terms, 
the so-called generalized Korteweg-de Vries-Benjamin-Bona-Mahony-Burgers equation 
with a fourth-order dissipative term (see \cite{ras-nik-kho})
\begin{eqnarray}
 \left\{\begin{array}{ll}
  \partial_tu +\partial_x \big( \, 
f(u) \, \big) - \alpha \, \partial_t\partial_x^2 u 
- \beta \, \partial_x^2 u 
+\delta  \, \partial_x^3 u
+\gamma  \, \partial_x^4 u=0
  \qquad &\big( \, t>0, \: x\in \mathbb{R} \, \big), \\[5pt]
 u(0, x) = u_0(x) \rightarrow u_{\pm} \qquad &( x \rightarrow \pm \infty ), 
 \end{array}
 \right.\,
\end{eqnarray}
where $\delta \in \mathbb{R}$, we can obtain the same stabilities as Theorems 1.1 and 1.2 
in the next theorems. 

\medskip

\noindent
{\bf Theorem 1.3} (Main Theorem I\hspace{-.1em}I\hspace{-.1em}I){\bf .}\quad{\it
Assume the far field states $u_{\pm}$ satisfy $u_- = u_+=\tilde{u}$, and 
the convective flux $f \in C^2(\mathbb{R})$. 
Further assume the initial data satisfy
$u_0-\tilde{u} \in L^2$ and
$\partial _xu_0 \in H^2$. 
Then for the Cauchy problem (1.9), 
the same result as in Theorem 1.1 holds true. 
}

\medskip

\noindent
{\bf Theorem 1.4} (Main Theorem I\hspace{-.1em}V){\bf .}\quad{\it
Assume the far field states $u_{\pm}$ satisfy $u_- < u_+$, and 
the convective flux $f \in C^5(\mathbb{R})$ 
satisfy {\rm(1.6)}. 
Further assume the initial data satisfy
$u_0-u_0 ^{\rm{R}} \in L^2$ and
$\partial _xu_0 \in H^2$. 
Then for the Cauchy problem (1.9), 
the same result as in Theorem 1.2 holds true. 
}

\medskip

We finally remark that when $\gamma=\delta=0$, 
the problem (1.9) becomes the following Cauchy problem for 
usual generalized Benjamin-Bona-Mahony-Burgers equation. 
\begin{eqnarray}
 \left\{\begin{array}{ll}
  \partial_tu +\partial_x \big( \, 
f(u) \, \big) - \alpha \, \partial_t\partial_x^2 u 
- \beta \, \partial_x^2 u =0
  \qquad &\big( \, t>0, \: x\in \mathbb{R} \, \big), \\[5pt]
 u(0, x) = u_0(x) \rightarrow u_{\pm} \qquad &( x \rightarrow \pm \infty ), 
 \end{array}
 \right.\,
\end{eqnarray}

\medskip

We can also obtain the similar stabilities as Theorems 1.1 and 1.2 
in the next theorems. 

\medskip

\noindent
{\bf Theorem 1.5.}\quad{\it
Assume the far field states $u_{\pm}$ satisfy $u_- = u_+=\tilde{u}$, and 
the convective flux $f \in C^2(\mathbb{R})$. 
Further assume the initial data satisfy
$u_0-\tilde{u} \in L^2$ and
$\partial _xu_0 \in H^2$. 
Then the Cauchy problem {\rm(1.10)} has a 
unique global in time 
solution $u$ 
satisfying 
\begin{eqnarray*}
 \left\{\begin{array}{ll}
u-\tilde{u} \in C^0\bigl( \, [\, 0, \, \infty \, ) \, ; H^3 \, \bigr), \\[5pt]
\partial _x u \in L^2\bigl( \, 0,\, \infty \, ; H^2 \, \bigr), \\[5pt]
\partial _t u \in L^2\bigl( \, 0,\, \infty \, ; H^2 \, \bigr),
\end{array} 
\right.\,
\end{eqnarray*}
and the asymptotic behavior 
$$
\lim _{t \to \infty}
\bigg( \,
\sup_{x\in \mathbb{R}}
| \, u(t,x) - \tilde{u} \, | 
+\sup_{x\in \mathbb{R}}
 |\,\partial_xu (t, x)\,| 
+\sup_{x\in \mathbb{R}}
 |\,\partial_x^2u (t, x)\,|
\, \bigg)= 0.
$$
}

\medskip

\noindent
{\bf Theorem 1.6.}\quad{\it
Assume the far field states $u_{\pm}$ satisfy $u_- < u_+$, and 
the convective flux $f \in C^5(\mathbb{R})$ 
satisfy {\rm(1.6)}. 
Further assume the initial data satisfy
$u_0-u_0 ^{\rm{R}} \in L^2$ and
$\partial _xu_0 \in H^2$. 
Then the Cauchy problem {\rm(1.10)} has a 
unique global in time 
solution $u$ 
satisfying 
\begin{eqnarray*}
 \left\{\begin{array}{ll}
u-u_0 ^{\rm{R}} \in C^0\bigl( \, [\, 0, \, \infty \, ) \, ; H^3 \, \bigr), \\[5pt]
\partial _x u \in L^2_{\rm{loc}}\bigl( \, 0,\, \infty \, ; H^2 \, \bigr), \\[5pt]
\partial _t u \in L^2_{\rm{loc}}\bigl( \, 0,\, \infty \, ; H^2 \, \bigr),
\end{array} 
\right.\,
\end{eqnarray*}
and the asymptotic behavior 
\begin{eqnarray*}
\left\{\begin{array}{ll}
\displaystyle{
\lim _{t \to \infty}\sup_{x\in \mathbb{R}} \, 
\left| \, 
u(t,x) - u^r\left(\, \frac{x}{t}\: ;\: u_-, \, u_+ \, \right) 
\, \right|} = 0,\\[10pt]
\displaystyle{
\lim _{t \to \infty}\sup_{x\in \mathbb{R}} \, 
| \, 
\partial_{x}u(t,x) 
- \partial_{x}u^r\left(\, t, \, x\: ;\: u_-, \, u_+ \, \right)
\, |} = 0, \\[10pt]
\displaystyle{
\lim _{t \to \infty}\sup_{x\in \mathbb{R}} \, 
| \, 
\partial_{x}^2u(t,x) 
- \partial_{x}^2u^r\left(\, t, \, x\: ;\: u_-, \, u_+ \, \right)
\, |} = 0, 
\end{array} 
\right.\,
\end{eqnarray*}
where $\partial_{x}u^r$ and $\partial_{x}^2u^r$ 
are given in Theorem 1.2.
}

Because the proofs of Theorems 1.1 and 1.3-1.6 
are similarly given as or easier than that for Theorem 1.2, 
we only show Theorem 1.2 in the following sections. 

\medskip

This paper is organized as follows. 
In Section 2, we construct the approximation of the rarefaction wave
and prepare the basic
properties of the rarefaction wave and the approximated one. 
We reformulate the problem 
in terms of the deviation from 
the asymptotic state in Section 3.
In order to show the asymptotics, 
we establish the {\it a priori} estimates 
by using the technical energy method in Section 4.  
Finally in Section 5, we give several uniform estimates 
by using the {\it a priori} estimates in Sections 3 and 4.

\medskip

{\bf Some Notation.}\quad 
We denote by $C$ generic positive constants unless 
they need to be distinguished. 
In particular, use 
$C_{\alpha, \beta, \cdots }$ 
when we emphasize the dependency on $\alpha,\: \beta,\: \cdots $.

For function spaces, 
${L}^p = {L}^p(\mathbb{R})$ and ${H}^k = {H}^k(\mathbb{R})$ 
denote the usual Lebesgue space and 
$k$-th order Sobolev space on the whole space $\mathbb{R}$ 
with norms $||\cdot||_{{L}^p}$ and $||\cdot||_{{H}^k}$, 
respectively. 

\bigskip 

\noindent
\section{Preliminaries} 
In this section, 
we 
prepare the several lemmas concerning with 
the basic properties of 
the rarefaction wave 
for the proof of the main Theorem 1.2. 
Since the rarefaction wave $u^r$ is not smooth enough, 
we construct a smooth approximated one. 
To do that, we first consider the rarefaction wave solution $w^r$ 
to the Riemann problem 
for the non-viscous Burgers equation 
\begin{equation}
\label{riemann-burgers}
  \left\{\begin{array}{l}
  \partial _t w + 
  \displaystyle{ \partial _x \left( \, \frac{1}{2} \, w^2 \right) } = 0 
  \, \, \; \; \qquad \quad \qquad \big( \,  t > 0,\: x\in \mathbb{R} \, \big),\\[7pt]
  w(0, x) = w_0 ^{\rm{R}} ( \, x\: ;\: w_- ,\: w_+\, ):= \left\{\begin{array}{ll}
                                                  w_+ & \, \: \; \quad (x>0),\\[5pt]
                                                  w_- & \, \: \; \quad (x<0),
                                                  \end{array}
                                                  \right.
  \end{array}
  \right.
\end{equation}
where $w_\pm \in \mathbb{R} \: (w_-<w_+)$ are 
the prescribed far field states. 
The unique global weak solution 
$w=w^r\left( \, x/t\: ;\: w_-,\: w_+ \, \right)$ 
of (\ref{riemann-burgers}) is explicitly given by 
\begin{equation}
\label{rarefaction-burgers}
w^r \left( \, \frac{x}{t}\: ;\: w_-,\: w_+ \, \right) 
= 
  \left\{\begin{array}{ll}
  w_{-} & \bigl(\, x \leq w_{-} \, t \, \bigr),\\[5pt]
  \displaystyle{ \frac{x}{t} } & \bigl(\, w_{-}\, t \leq x \leq w_{+}\, t \, \bigr),\\[5pt]
  w_+ & \bigl(\, x\geq w_{+}\, t \, \bigr).
  \end{array}\right.
\end{equation} 
Next, under the condition 
$f''(u)>0\ (u\in \mathbb{R})$ and $u_-<u_+$, 
the rarefaction wave solution 
$u=u^r\left( \, \frac{x}{t}\: ;\: u_-,\: u_+ \, \right)$ 
of the Riemann problem (1.2) 
for hyperbolic conservation law 
is exactly given by 
\begin{equation}
u^r\left( \, \frac{x}{t} \: ; \:  u_-,\: u_+ \, \right) 
= (\lambda)^{-1}
  \biggl( \, w^r
  \left( \, 
  \frac{x}{t} \: ; \:  \lambda_-,\: \lambda_+ \, \right) \, \biggr)
\end{equation}
which is nothing but (1.6), 
where $\lambda_\pm := \lambda(u_\pm) = f'(u_\pm)$. 
We define a smooth approximation of 
$w^r(\, \frac{x}{t}\: ;\: w_-,\: w_+ \,)$ 
by the unique classical solution 
$$
w=w(\, t, \, x\: ;\: w_-,\: w_+ \,)
$$
to the Cauchy problem for the following 
non-viscous Burgers equation as 
\begin{eqnarray}
\label{smoothappm}
\left\{\begin{array}{l}
 \partial _t w 
 + \displaystyle{ \partial _x \left( \, \frac{1}{2} \, w^2 \, \right) } =0 
 \, \, \; \; \quad \qquad \qquad \qquad \qquad \qquad  \qquad  \quad  \:
 \big( \, t>0,\: x\in \mathbb{R} \, \big),\\[7pt]
 w(0,x) 
 = w_0(x) 
 := \displaystyle{ \frac{w_-+w_+}{2} 
    + \frac{w_+-w_-}{2}\,K_{q} \, \int_{0}^{\epsilon x} \frac{\mathrm{d}y}{(1+y^2)^q} }
 \quad \quad \; (x\in \mathbb{R}), 
\end{array}
\right.
\end{eqnarray}   
where $K_{q}$ is a positive constant such that 
$$
K_{q} \, \int_{0}^{\infty} \frac{\mathrm{d}y}{(1+y^2)^q} =1 
\quad \left( \, q>\frac{1}{2} \, \right). 
$$
By applying the method of characteristics, 
we get the following formula 
\begin{eqnarray}
 \left\{\begin{array} {l}
 w(t, x)=w_0\bigl( \, x_0(t, x) \, \bigr)=
 \displaystyle{ \frac{\lambda_-+\lambda_+}{2} } 
+ \displaystyle{ \frac{\lambda_+-\lambda_-}{2}
\,K_{q} \, \int_{0}^{\epsilon x_0(t, x) } \frac{\mathrm{d}y}{(1+y^2)^q}  } ,\\[7pt]
 x=x_0(t, x)+w_0\bigl( \, x_0(t, x) \, \bigr)\,t.
 \end{array}
  \right.\,
\end{eqnarray}
By making use of (2.5) similarly as in \cite{matsu-nishi1}, 
we can obtain the properties of 
the smooth approximation $w(\, t, \, x\: ;\: w_-, \: w_+ \,)$ 
in the next lemma.

\medskip

\noindent
{\bf Lemma 2.1.}\quad{\it
Assume that the far field states satisfy $w_-<w_+$. 
Then the classical solution $w(t,\, x)=w(\, t, \, x\: ;\: w_-,\: w_+ \,)$
given by {\rm(2.4)} 
satisfies the following properties: 

\noindent
{\rm (1)}\ \ $w_- < w(t, x) < w_+$ and\ \ $\partial_xw(t, x) > 0$  
\quad  $\big( \, t>0, \: x\in \mathbb{R} \, \big)$.

\smallskip

\noindent
{\rm (2)}\ For any $r \in [\, 1, \, \infty \,]$, there exists a positive 
constant $C_{q, r}$ such that
             \begin{eqnarray*}
                 \begin{array}{l}
                    \| \, \partial_x w(t) \, \|_{L^r}^r \leq 
                    C_{q,r} \, 
                    \min \left\{ \, 
                    \epsilon^{r-1} \, 
                    \tilde{w}^{r}, \, 
                    \tilde{w}\, 
                    ( 1+t )^{-r+1}
                    \, \right\}
                    \quad \bigl( \, t\ge 0 \, \bigr),\\[5pt]
                    \| \, \partial_t w(t) \, \|_{L^r}^r \leq 
                    C_{q,r} \, \tilde{\tilde{w}}^{r}\, 
                    \min \left\{ \, 
                    \epsilon^{r-1} \, 
                    \tilde{w}^{r}, \, 
                    \tilde{w}\, 
                    ( 1+t )^{-r+1}
                    \, \right\}
                    \quad \bigl( \, t\ge 0 \, \bigr),\\[5pt]
                    \| \, \partial_x^2 w(t) \, \|_{L^r}^r \\[5pt]
                    \leq 
                    C_{q,r} \, 
                    \min \left\{ \, 
                    \epsilon^{2r-1} \, 
                    \tilde{w}^{r}, \, 
                    \epsilon^{(r-1)(1-\frac{1}{2q})} \, 
                    \tilde{w}^{-\frac{r-1}{2q}}\, 
                    ( 1+t )^{-r-\frac{r-1}{2q}}
                    \, \right\}
                    \quad \bigl( \, t\ge 0 \, \bigr),\\[5pt]
                    \| \, \partial_x^3 w(t) \, \|_{L^r}^r \leq 
                    C_{q,r} \, 
                    \min \left\{ \, 
                    \epsilon^{3r-1} \, 
                    \tilde{w}^{r}, \, 
                    a(1+t, \epsilon, \tilde{w})
                    \, \right\} 
                    \quad \bigl( \, t\ge 0 \, \bigr),\\[5pt]
                    \| \, \partial_x^4 w(t) \, \|_{L^r}^r \leq 
                    C_{q,r} \, 
                    \min \left\{ \, 
                    \epsilon^{4r-1} \, 
                    \tilde{w}^{r}, \, 
                    b(1+t, \epsilon, \tilde{w})
                    \, \right\} 
                    \quad \bigl( \, t\ge 0 \, \bigr),
                    \end{array}       
              \end{eqnarray*}
where 
\begin{equation*}
\tilde{w}:= \frac{w_{+}-w_{-}}{2}>0,
\quad \tilde{\tilde{w}}:= \max \{ \, | \, w_{-} \, |, \, | \, w_{+} \, | \, \}, 
\end{equation*}
\begin{align*}
\begin{aligned}
a(t, \epsilon, \tilde{w})
&:=\epsilon^{3r} \, 
                    \tilde{w}^{r}\, 
                    ( 1+\epsilon \, \tilde{w} \, t )^{1-4r}
  +\epsilon^{2(r-1)(1-\frac{1}{2q})} \, 
                    \tilde{w}^{-\frac{r-1}{q}}\, 
                    t^{-1-(r-1)\left( 1+ \frac{1}{q}\right)}\\
 &\quad 
  +\epsilon^{(2r-1)(1-\frac{1}{2q})} \, 
                    \tilde{w}^{-\frac{2r-1}{2q}}\, 
                    t^{-r-\frac{2r-1}{2q}},
\end{aligned}      
\end{align*}
and 
\begin{align*}
\begin{aligned}
b(t, \epsilon, \tilde{w})
&:=\epsilon^{3r} \, 
                    \tilde{w}^{r}\, 
                    ( 1+\epsilon \, \tilde{w} \, t )^{1-5r}
  +\epsilon^{(3r-2)(1-\frac{1}{2q})} \, 
                    \tilde{w}^{-\frac{3r-2}{2q}}\, 
                    t^{-r \left( 1+ \frac{3}{2q}\right)+ \frac{1}{q}}\\
 &\quad 
  +\epsilon^{(3r-1)(1-\frac{1}{2q})} \, 
                    \tilde{w}^{-\frac{3r-1}{2q}}\, 
                    t^{-r-\frac{3r-1}{2q}}.
\end{aligned}      
\end{align*}
              
\smallskip

\noindent
{\rm (3)}\ It follows that 
\begin{eqnarray*}
\left\{\begin{array}{ll}
\displaystyle{
\lim _{t \to \infty}\sup_{x\in \mathbb{R}} \, 
\left| \, w(t, x)- w^r \left( \frac{x}{t} \right) 
\, \right|} = 0,\\[10pt]
\displaystyle{
\lim _{t \to \infty}\sup_{x\in \mathbb{R}} \, 
| \, 
\partial_{x}w(t,x) 
- \partial_{x}w^r(t, x)
\, |} = 0,
\end{array} 
\right.\,
\end{eqnarray*}
where $\partial_{x}w^r$ 
is given by 
\begin{equation*}
\partial_{x}w^r\left(\, t, \, x \: ;\: w_-, \, w_+ \, \right)
= \left\{
\begin{array}{ll}
  0  & \; \bigl(\, x \leq w_-\,t \, \bigr),\\[7pt]
  \displaystyle{\frac{1}{t}}
  & \; \bigl(\, w_-\,t \leq x \leq w_+\,t\,  \bigr),\\[15pt]
  0 & \; \bigl(\, x \geq w_+\,t \, \bigr).
\end{array}
\right. 
\end{equation*}
}

\bigskip

We now define the approximation for 
the rarefaction wave $u^r\left( \, {x}/{t}\: ;\: u_-,\: u_+ \, \right)$ by 
\begin{equation}
U^r(\, t, \, x\: ; \: u_-,\: u_+ \,) 
= (\lambda)^{-1} \bigl( \,w(\, t, \, x\: ;\: \lambda_-,\: \lambda_+\,) \, \bigr).
\end{equation}
Using Lemma 2.1, 
we also have the next lemma. 

\medskip

\noindent
{\bf Lemma 2.2.}\quad{\it
Let $q=1$. Assume that the far field states satisfy $u_-<u_+$, 
and the flux function $f\in C^5(\mathbb{R})$, 
$f''(u)>0 \; (\,u\in [\,u_-,\,u_+\,]\,)$. 
Then we have the following properties. 

\noindent
{\rm (1)}\ $U^r(t, x)$ defined by {\rm (2.6)} is 
the unique $C^4$-global solution in space-time 
of the Cauchy problem
$$
\left\{
\begin{array}{l} 
\partial _t U^r +\partial _x \bigl( \, f(U^r ) \, \bigr) = 0 
\quad 
\big( \, t>0, \: x\in \mathbb{R} \, \big),\\[7pt]
U^r(0,x) 
= \displaystyle{ (\lambda)^{-1} \left( \, \frac{\lambda_- + \lambda_+}{2} 
+ \frac{\lambda_+ - \lambda_-}{2} 
  \,K_{q} \, \int_{0}^{\epsilon x} \frac{\mathrm{d}y}{(1+y^2)^q} \, \right) } 
\quad( x\in \mathbb{R}),\\[10pt]
\displaystyle{\lim_{x\to \pm \infty}} U^r(t, x) =u_{\pm} 
\quad 
\bigl( \, t\ge 0 \, \bigr).
\end{array}
\right.\,     
$$
{\rm (2)}\ \ $u_- < U^r(t, x) < u_+$ 
and\ \ $\partial_xU^r(t,x) > 0$  
\quad  $\bigl( \, t>0, \: x\in \mathbb{R} \, \bigr)$.

\smallskip

\noindent
{\rm (3)}\ For any $r \in [\,1,\, \infty \,]$, there exists a positive 
constant $C_{\lambda_{\pm},r}$ such that
             \begin{eqnarray*}
                 \begin{array}{l}
                    \| \, \partial_x U^r(t) \, \|_{L^r}^r 
                    \leq C_{\lambda_{\pm},r}\, 
                    \min \left\{ \, 
                    \epsilon^{r-1}, \, 
                    ( 1+t )^{-r+1}
                    \, \right\} 
                    \quad \bigl(\, t\ge 0 \, \bigr),\\[5pt]
                    \| \, \partial_t U^r(t) \, \|_{L^r}^r 
                    \leq C_{\lambda_{\pm},r}\, 
                    \min \left\{ \, 
                    \epsilon^{r-1}, \, 
                    ( 1+t )^{-r+1}
                    \, \right\} 
                    \quad \bigl(\, t\ge 0 \, \bigr),\\[5pt]
                    \| \, \partial_x^2 U^r(t) \, \|_{L^r}^r 
                    \leq C_{\lambda_{\pm},r}\, 
                    \min \left\{ \, 
                    \epsilon^{2r-1}, \, 
                    \epsilon^{\frac{r-1}{2}} \, 
                    ( 1+t )^{-\frac{3r-1}{2}}
                    \, \right\}
                    \quad \bigl(\, t\ge 0 \, \bigr),\\[5pt]
                    \| \, \partial_x^3 U^r(t) \, \|_{L^r}^r 
                    \leq C_{\lambda_{\pm},r}\, 
                    \min \left\{ \, 
                    \epsilon^{3r-1}, \, 
                    \epsilon^{r-1} \, 
                    ( 1+t )^{-2r+1}
                    \, \right\}
                    \quad \bigl(\, t\ge 0 \, \bigr),\\[5pt]
                    \| \, \partial_x^4 U^r(t) \, \|_{L^r}^r 
                    \leq C_{\lambda_{\pm},r}\, 
                    \min \left\{ \, 
                    \epsilon^{4r-1}, \, 
                    \epsilon^{\frac{3r-2}{2}} \, 
                    ( 1+t )^{-\min \left\{ \frac{3}{2} r +1 , 
                            \,  \frac{5}{2} r -1 \right\}}
                    \, \right\}
                    \quad \bigl(\, t\ge 0 \, \bigr).
                    \end{array}       
              \end{eqnarray*}
              
\smallskip

\noindent
{\rm (4)}\ It follows that 
\begin{eqnarray*}
\left\{\begin{array}{ll}
\displaystyle{\lim_{t\to \infty} 
\sup_{x\in \mathbb{R}}
\left| \,U^r(t, x)- u^r \left( \frac{x}{t} \right) \, \right| = 0},\\[10pt]
\displaystyle{
\lim _{t \to \infty}\sup_{x\in \mathbb{R}} \, 
| \, 
\partial_{x}U^r(t,x) 
- \partial_{x}u^r(t, x)
\, |} = 0,\\[10pt]
\displaystyle{
\lim _{t \to \infty}\sup_{x\in \mathbb{R}} \, 
\left| \, 
\partial_{x}^2U^r(t,x) 
- \partial_{x}^2u^r(t, x)
\, \right|} = 0.
\end{array} 
\right.\,
\end{eqnarray*}
%
%
%
%
%
%
%
}

\medskip

Because the proofs of Lemmas 2.1 and 2.2 are given in 
\cite{hashimoto-matsumura}, \cite{hattori-nishihara}, \cite{liu-matsumura-nishihara}, \cite{matsu-nishi1}, 
\cite{matsumura-yoshida}, \cite{wan-zhu}, \cite{yoshida1}, \cite{yoshida7}, \cite{yoshida10} and so on, 
we omit the proofs here. 

%

\bigskip 

\noindent
\section{Reformulation of the problem} 
In this section, we reformulate our problem (1.1) 
in terms of the deviation from the asymptotic state. 
Putting $\phi$ as 
\begin{equation}
u(t, x) = U^r(t, x) + \phi(t, x), 
\end{equation}
we reformulate the problem (1.1) in terms of 
the deviation $\phi $ from $U^r$ as 
\begin{eqnarray}
 \left\{\begin{array}{ll}
  \partial _t\phi 
  + \partial_x \big( \, f(\phi+U^r) - f(U^r) \, \big) \\[5pt]
  \quad 
  - \alpha \, \partial_t\partial_x^2 \phi 
  - \beta \, \partial_x^2 \phi 
  + \gamma \, \partial_x^4 \phi 
    = F(U^r)
     \quad  \bigl( \, t>0,\: x\in \mathbb{R} \, \bigr), \\[5pt]
  \phi(0, x) = \phi_0(x) 
  := u_0(x)-U^r(0, x) %
  \rightarrow 0 \quad (x\rightarrow \pm \infty),
 \end{array}
 \right.\,
\end{eqnarray}
where 
$$
F(U^r):= \alpha \, \partial_t\partial_x^2 U^r
  +\beta \, \partial_x^2 U^r 
  - \delta \, \partial_x^3 U^r - \gamma \, \partial_x^4 U^r. 
$$
Then we look for 
the unique global in time solution 
$\phi $ which has the asymptotic behavior 
\begin{equation}
\displaystyle{
\sum_{k=1}^3 \, 
\sup_{x\in \mathbb{R}}
 \left|\,\partial_x^k\phi (t, x)\, \right| \to  0 
\qquad (t\to  \infty)}. 
\end{equation}
Here we note 
that $\phi_0 \in H^3$ by the assumptions on $u_0$ 
and Lemma 2.2. 
Then the corresponding theorems 
for $\phi$ to Theorem 1.2 we should prove is 
as follows. 

\medskip

\noindent
{\bf Theorem 3.1.} (Global Existence){\bf .}\quad{\it
Assume the far field states $u_{\pm}$ satisfy $u_- < u_+$, and 
the convective flux $f \in C^5(\mathbb{R})$ 
satisfy {\rm(1.6)}. 
Further assume the initial data satisfy
$\phi_0 \in H^3$. 
Then the Cauchy problem {\rm(3.2)} has a 
unique global in time 
solution $\phi$ 
satisfying 
\begin{eqnarray*}
\left\{\begin{array}{ll}
\phi \in C^0\bigl( \, [\, 0, \, \infty ) \, ; H^3 \, \bigr), \\[5pt]
\partial _x \phi \in L^2\bigl( \, 0,\, \infty \, ; H^3 \, \bigr), \\[5pt]
\partial _t \phi \in L^2\bigl( \, 0,\, \infty \, ; H^2 \, \bigr),
\end{array} 
\right.\,
\end{eqnarray*}
and the asymptotic behavior 
$$
\lim _{t \to \infty} \,
\sum_{k=0}^2 \, 
\sup_{x\in \mathbb{R}}
 \left|\,\partial_x^k\phi (t, x)\, \right|= 0. 
$$
}

\medskip

%
%
In order to obtain Theorem 3.1, 
we prepare the local existence precisely, 
we formulate 
the problem (3.2) at general 
initial time $\tau \ge 0$: 
\begin{eqnarray}
 \left\{\begin{array}{ll}
  \partial _t\phi 
  + \partial_x \big( \, f(\phi+U^r) - f(U^r) \, \big) \\[5pt]
  \quad 
  - \alpha \, \partial_t\partial_x^2 \phi 
  - \beta \, \partial_x^2 \phi 
  + \gamma \, \partial_x^4 \phi 
    = F(U^r)
     \quad  \bigl( \, t>\tau,\: x\in \mathbb{R} \, \bigr), \\[5pt]
  \phi(\tau, x) = \phi_\tau(x) 
  := u_\tau(x)-U^r(\tau, x) %
  \rightarrow 0 \quad (x\rightarrow \pm \infty).
 \end{array}
 \right.\,
\end{eqnarray}

\medskip

\noindent
{\bf Theorem 3.2} (Local Existence){\bf .}\quad{\it
For any $ M > 0 $, there exists a positive constant 
$t_0=t_0(M)$ not depending on $\tau$ 
such that if 
$\phi_{\tau} \in H^3$ and 
$$
\displaystyle{ 
\| \, \phi_{\tau} \, \|_{H^3} \leq M}, 
$$
then the Cauchy problem {\rm (3.4)} has a unique solution $\phi$ 
on the time interval $[\, \tau, \, \tau+t_{0}(M)\, ]$ satisfying 
\begin{eqnarray*}
\left\{\begin{array}{ll}
\phi \in C^0
\bigl( \, [\, \tau, \, \tau+t_{0}\, ] \, ; H^3 \bigr),\\[5pt]
\partial _x \phi  \in L^2\bigl( \, \tau, \, \tau+t_{0} \, ; H^3 \bigr),\\[5pt]
\partial _t \phi  \in L^2\bigl( \, \tau, \, \tau+t_{0} \, ; H^2 \bigr),\\[5pt]
\displaystyle{
\sup_{t \in [\tau, \tau + t_{0}]} \, 
\| \, \phi(t) \, \|_{H^3}
}
\leq 2 \, M. 
\end{array} 
\right.\,
\end{eqnarray*}
}

\medskip

\noindent
Because the proof of Theorem 3.2 is standard, 
we omit the details here (cf. \cite{yin-zhao-kim}, \cite{zhao-xuan}). 
The {\it a priori} estimates we establish in Section 4 
are the following. 

\medskip

\noindent
{\bf Theorem 3.3} ({\it A Priori} Estimates){\bf .}\quad{\it 
Under the same assumptions as in Theorem 3.1, 
for any initial data 
$\phi_0 \in H^3$, 
there exists a positive constant $C_{\phi_0}$
such that 
if the Cauchy problem {\rm (3.1)}
has a solution 
$\phi$ 
on the time interval $[\, 0, \, T\, ]$ satisfying 
\begin{eqnarray*}
\left\{\begin{array}{ll}
\phi \in C^0
\bigl( \, [\, 0, \, T\, ] \, ; H^3 \bigr),\\[5pt]
\partial _x \phi  \in L^2\bigl( \, 0, \, T \, ; H^3 \bigr),\\[5pt]
\partial _t \phi  \in L^2\bigl( \, 0, \, T \, ; H^2 \bigr),
\end{array} 
\right.\,
\end{eqnarray*}
for some positive constant $T$, 
then it holds that 
\begin{align}
\begin{aligned}
& \| \, \phi(t) \, \|_{H^3}^{2} 
+ \int^{t}_{0} 
   \big\| \, 
   \big(\, \sqrt{\partial_x U^r} \:  \phi \, \big)(\tau) 
   \, \big\|_{L^2}^{2} \, \mathrm{d}\tau 
  \\
& 
+  \int^{t}_{0} 
\| \, \partial_{x}\phi(\tau) \, \|_{H^3}^{2} \, \mathrm{d}\tau
+ \int^{t}_{0} 
\| \, \partial_{t}\phi(\tau) \, \|_{H^2}^{2} \, \mathrm{d}\tau  
\leq C_{\phi_0} 
\quad \big( \, t \in [\, 0, \, T\, ] \, \big).
\end{aligned}
\end{align}
}

\medskip

Combining the local existence Theorem 3.2 together with 
the {\it a priori} estimates, Theorem 3.3, we can obtain 
global existence Theorem 3.1. 
In fact, we can obtain the unique global in time solutions $\phi$ to (3.2) in 
Theorem 3.1 satisfying 
\begin{eqnarray*}
\left\{\begin{array}{ll}
\phi \in C^0\bigl( \, [\, 0, \, \infty ) \, ; H^3 \, \bigr), \\[5pt]
\partial _x \phi \in L^2\bigl( \, 0,\, \infty \, ; H^3 \, \bigr), \\[5pt]
\partial _t \phi \in L^2\bigl( \, 0,\, \infty \, ; H^2 \, \bigr),
\end{array} 
\right.\,
\end{eqnarray*}
and 
\begin{align}
\begin{aligned}
& \sup_{t \ge 0}
  \| \, \phi(t) \, \|_{H^3}^{2} 
 + \int^{\infty}_{0} 
   \big\| \, 
   \big(\, \sqrt{\partial_x U^r} \:  \phi \, \big)(t) 
   \, \big\|_{L^2}^{2} \, \mathrm{d}t
  \\
& \quad
 + \int^{\infty}_{0} 
\| \, \partial_{x}\phi(t) \, \|_{H^3}^{2} \, \mathrm{d}t
+ \int^{\infty}_{0} 
\| \, \partial_{t}\phi(t) \, \|_{H^2}^{2} \, \mathrm{d}t
< \infty 
\end{aligned}
\end{align}
which yields 
\begin{equation}
\int _0^{\infty }\bigg|\,\frac{\mathrm{d}}{\mathrm{d}t} 
\|  \, \partial_x\phi(t)  \, \|_{L^2}^2 \,\bigg| \, \mathrm{d}t
< \infty, \quad 
\int _0^{\infty }\bigg|\,\frac{\mathrm{d}}{\mathrm{d}t} 
\|  \, \partial_x^2\phi(t)  \, \|_{L^2}^2 \,\bigg| \, \mathrm{d}t
< \infty.
\end{equation}
%
%
%
%
In fact, by using (3.6), direct computation shows that 
\begin{align}
\begin{aligned}
\int _0^{\infty }\bigg|\,\frac{\mathrm{d}}{\mathrm{d}t} 
\|  \, \partial_x\phi(t)  \, \|_{L^2}^2 \,\bigg| \, \mathrm{d}t
&= 2\, \int _0^{\infty } \left| \, \int _{-\infty}^{\infty } 
\partial_x\phi \, \partial_t\partial_x\phi 
\, \mathrm{d}x \, \right|
    \, \mathrm{d}t\\
& \le \int _0^{\infty } \left( \, 
      \|  \, \partial_x\phi(t)  \, \|_{L^2}^2
      +\|  \, \partial_t\partial_x\phi(t)  \, \|_{L^2}^2
      \, \right)  \, \mathrm{d}t< \infty, 
\end{aligned}
\end{align}
\begin{align}
\begin{aligned}
\int _0^{\infty }\bigg|\,\frac{\mathrm{d}}{\mathrm{d}t} 
\|  \, \partial_x^2\phi(t)  \, \|_{L^2}^2 \,\bigg| \, \mathrm{d}t
&= 2\, \int _0^{\infty } \left| \, \int _{-\infty}^{\infty } 
\partial_x^2\phi \, \partial_t\partial_x^2\phi 
\, \mathrm{d}x \, \right|
    \, \mathrm{d}t\\
& \le \int _0^{\infty } \left( \, 
      \|  \, \partial_x^2\phi(t)  \, \|_{L^2}^2
      +\|  \, \partial_t\partial_x^2\phi(t)  \, \|_{L^2}^2
      \, \right)  \, \mathrm{d}t< \infty. 
\end{aligned}
\end{align}
From (3.8) and (3.9), we get (3.7).
We immediately have from (3.7) that 
\begin{equation}
\lim_{t \to \infty} \|  \, \partial_x\phi(t)  \, \|_{L^2} = 0, \quad 
\lim_{t \to \infty} \|  \, \partial_x^2\phi(t)  \, \|_{L^2} = 0.
\end{equation}
Further from (3.10), by using the Sobolev inequality, 
we obtain the desired asymptotic behavior (3.3) as follows. 
\begin{equation*}
\sup_{x \in \mathbb{R}}|\,\phi(t,x)\,|
\leq \sqrt{2}\, \| \,\phi(t)\, \|^{\frac{1}{2}}_{L^2} 
\| \,\partial_x\phi(t)\, \| ^{\frac{1}{2}}_{L^2}\rightarrow 0
\quad \big( \, t \rightarrow \infty \, \big), 
\end{equation*} 
\begin{equation*}
\sup_{x \in \mathbb{R}}|\,\partial_x\phi(t,x)\,|
\leq \sqrt{2}\, \| \,\partial_x\phi(t)\, \|^{\frac{1}{2}}_{L^2} 
\| \,\partial_x^2\phi(t)\, \| ^{\frac{1}{2}}_{L^2}\rightarrow 0
\quad \big( \, t \rightarrow \infty \, \big), 
\end{equation*} 
\begin{equation*}
\sup_{x \in \mathbb{R}}|\,\partial_x^2\phi(t,x)\,|
\leq \sqrt{2}\, \| \,\partial_x^2\phi(t)\, \|^{\frac{1}{2}}_{L^2} 
\| \,\partial_x^3\phi(t)\, \| ^{\frac{1}{2}}_{L^2}\rightarrow 0
\quad \big( \, t \rightarrow \infty \, \big). 
\end{equation*} 

Thus Theorems 3.1 is proved.


\bigskip 

\noindent
\section{{\it A priori} estimates}
In this section, 
we show the following {\it a priori} estimate for $\phi$ in Theorem 3.3. 
To do that, we prepare the following basic estimate. 

\medskip

\noindent
{\bf Proposition 4.1.}\quad {\it
There exists a positive constant $C_{\phi_0}$
such that 
\begin{align*}
\begin{aligned}
& \| \, \phi(t) \, \|_{H^1}^{2} 
  +  \int^{t}_{0} \int^{\infty}_{-\infty} \int^{\phi}_{0} 
   \bigl( \, 
    f'(\eta +U^r)-f'(U^r) \,  \bigr) \, \mathrm{d}\eta 
    \, \partial_x U^r \, \mathrm{d}x \mathrm{d}\tau \\
& 
+ \int^{t}_{0} 
\| \, \partial_{x}\phi(\tau) \, \|_{H^1}^{2} \, \mathrm{d}\tau 
\leq C_{\phi_0} 
\quad \big( \, t \in [\, 0, \, T\, ] \, \big). 
\end{aligned}
\end{align*}
}

\medskip

{\bf Proof of Proposition 4.1.}
Multiplying the equation in (3.2) by $\phi$ and integrating it with respect to $x$, 
we have, after integration by parts, that 
\begin{align}
\begin{aligned}
& \frac{1}{2} \, \frac{\mathrm{d}}{\mathrm{d}t} \, \| \, \phi(t) \, \|_{L^2}^2 
+ \frac{\alpha}{2} \, \frac{\mathrm{d}}{\mathrm{d}t} \, \| \, \partial_x \phi(t) \, \|_{L^2}^2 \\
& + \int^{\infty}_{-\infty} \int^{\phi}_{0} 
   \bigl( \, 
    f'(\eta +U^r)-f'(U^r) \,  \bigr) \, \mathrm{d}\eta 
    \, \partial_x U^r \, \mathrm{d}x \\
&+ \beta \, \| \, \partial_{x}\phi(t) \, \|_{L^2}^{2} 
+ \gamma \, \| \, \partial_{x}^2\phi(t) \, \|_{L^2}^{2} 
= \int^{\infty}_{-\infty} \phi \, F(U^r) \, \mathrm{d}x.
\end{aligned}
\end{align}
By making use of the Sobolev inequality and the Young inequality, 
we estimate the right-hand-side of (4.1) as follows. 
\begin{align}
\begin{aligned}
\left| \,  \int^{\infty}_{-\infty} \phi \, F(U^r) \, \mathrm{d}x \, \right|
&\le \sqrt{2} \, \| \, \phi \, \|_{L^2}^{\frac{1}{2}}\, 
     \| \, \partial_x\phi \, \|_{L^2}^{\frac{1}{2}}\, \| \, F(U^r) \, \|_{L^1}\\
&\le \frac{\beta}{2} \, \| \, \partial_x\phi \, \|_{L^2}^2
     + C_{\beta} \, \| \, \phi \, \|_{L^2}^{\frac{2}{3}}
     \, \| \, F(U^r) \, \|_{L^1}^{\frac{4}{3}}\\
&\le \frac{\beta}{2} \, \| \, \partial_x\phi \, \|_{L^2}^2
     + C_{\beta} \, \left( \, 1+\| \, \phi \, \|_{L^2}^2 \, \right)
     \, \| \, F(U^r) \, \|_{L^1}^{\frac{4}{3}}.
\end{aligned}
\end{align}
Substituting (4.2) into (4.1), integrating the resultant inequality with respect to $t$, 
noting 
\begin{equation*}
\partial_t\partial_xU^r =-f''(U^r) \, |\, \partial_xU^r \, |^2-f'(U^r) \, \partial_x^2U^r,
\end{equation*}
\begin{equation*}
\| \, F(U^r) \, \|_{L^1}^{\frac{4}{3}} 
\le C_{\alpha,u_{\pm}} \, \| \, \partial_xU^r \, \|_{L^2}^{\frac{8}{3}}
+ C_{\alpha,\beta,u_{\pm}} \, \| \, \partial_x^2U^r \, \|_{L^1}^{\frac{4}{3}}
+\gamma^{\frac{4}{3}} \, \| \, \partial_x^4U^r \, \|_{L^1}^{\frac{4}{3}}
\in L_t^1(0, \, \infty)
\end{equation*}
from Lemma 2.2, 
and using the Gronwall inequality, we obtain the desired estimate. 

Thus, we complete the proof of Proposition 4.1. 

\medskip

Olny from Proposition 4.1, using the Sobolev inequality, 
we can easily get the uniform boundedness of $\phi$ in the next lemma
(cf. \cite{ilin-kalashnikov-oleinik}, \cite{kanel}, \cite{lad-sol-ura}). 

\medskip

\noindent
{\bf Lemma 4.2.}\quad {\it
There exists a positive constant $C_{\phi_0}$
such that 
\begin{equation*}
\sup_{t\in [0,T], \, x\in \mathbb{R}}\, 
| \, \phi (t,x) \, | \le C_{\phi_0}.
\end{equation*}
}

\medskip
By the uniform boundedness of $\phi$, Lemma 4.2, 
we note that the second term on the left-hand side of 
the {\it a priori} estimate in Proposition 4.1 can be replaced by 
the left-hand side of the following inequality as 
\begin{align}
\begin{aligned}
&\int^{t}_{0} \int^{\infty}_{-\infty} \int^{\phi}_{0} 
   \bigl( \, 
    f'(\eta +U^r)-f'(U^r) \,  \bigr) \, \mathrm{d}\eta 
    \, \partial_x U^r \, \mathrm{d}x \mathrm{d}\tau\\
&\geq C_{\phi_0}^{-1} \, \int^{t}_{0} \big\| \, 
   \big(\, \sqrt{\partial_x U^r} \:  \phi \, \big)(\tau) 
   \, \big\|_{L^2}^{2} \mathrm{d}\tau.
\end{aligned}
\end{align}

\medskip

We now prepare the following propieties for the {\it a priori} estimates 
for the derivatives of $\phi$
in the next lemma. 

\medskip

\noindent
{\bf Lemma 4.3.}\quad {\it
There exists a positive constant $C_{\phi_0}$
such that 
\begin{equation*}
\int_0^t 
\left\| \, \bigl( \, 
    f(\phi +U^r)-f(U^r) 
    \,  \bigr) (\tau)
    \, \right\|_{L^2}^2
    \, \mathrm{d}\tau
\leq C_{\phi_0} 
\quad \big( \, t \in [\, 0, \, T\, ] \, \big),
\end{equation*}
\begin{equation*}
\int_0^t 
\left\| \, F(U^r) (\tau)
    \, \right\|_{L^2}^2
    \, \mathrm{d}\tau
\leq C_{\phi_0} 
\quad \big( \, t \in [\, 0, \, T\, ] \, \big).
\end{equation*}
}

\medskip

{\bf Proof of Lemma 4.3.}
The second estimate, that is, the time-integrability of 
$\left\| \, F(U^r) \, \right\|_{L^2}^2$, is 
obtained by Lemma 2.2. 
Noting Proposition 4.1 and using Lagrange's mean-value theorem, 
we can also obtain the first estimate, that is, the time-integrability of 
$\left\| \, f(\phi +U^r)-f(U^r) \, \right\|_{L^2}^2$. 

Thus, the proof of Lemma 4.3 is completed. 

\medskip

\noindent
Next, we show the {\it a priori} 
estimate for $\partial_{x}\phi$ and $\partial_{x}^2\phi$ as follows. 

\medskip

\noindent
{\bf Proposition 4.4.}\quad {\it
There exists a positive constant $C_{\phi_0}$
such that 
\begin{equation*}
\| \, \partial_{x}\phi(t) \, \|_{H^1}^{2} 
  + \int^{t}_{0} 
  \left( \, 
  \| \, \partial_t\phi(\tau) \, \|_{H^1}^{2}
+\| \, \partial_{x}^2\phi(\tau) \, \|_{H^1}^{2} 
\, \right)
\, \mathrm{d}\tau
\leq C_{\phi_0} 
\quad \big( \, t \in [\, 0, \, T\, ] \, \big).
\end{equation*}
}

\medskip

{\bf Proof of Proposition 4.4.}
Multiplying the equation in (3.2) by $\partial_t\phi-\partial_x^2 \phi$ 
and integrating it with respect to $x$, 
we have, after integration by parts, that 
\begin{align}
\begin{aligned}
& \frac{\beta+1}{2} \, \frac{\mathrm{d}}{\mathrm{d}t} \, \| \, \partial_x \phi(t) \, \|_{L^2}^2 
+ \frac{\gamma+1}{2} \, \frac{\mathrm{d}}{\mathrm{d}t} \, \| \, \partial_x^2 \phi(t) \, \|_{L^2}^2 \\
&+ \| \, \partial_{t}\phi(t) \, \|_{L^2}^{2} 
+ \alpha \, \| \, \partial_{t}\partial_{x}\phi(t) \, \|_{L^2}^{2} 
+ \beta \, \| \, \partial_{x}^2\phi(t) \, \|_{L^2}^{2} 
+ \gamma \, \| \, \partial_{x}^3\phi(t) \, \|_{L^2}^{2} \\
&=-\int^{\infty}_{-\infty} 
\big( \, \partial_t\phi-\partial_x^2 \phi \, \big) \, 
   \bigl( \, 
    f(\phi +U^r)-f(U^r) 
    \,  \bigr) \, \mathrm{d}x \\
&+ \int^{\infty}_{-\infty} 
\big( \, \partial_t\phi-\partial_x^2 \phi \, \big) \, F(U^r) \, \mathrm{d}x.
\end{aligned}
\end{align}
By using the Young inequality, we estimate the each terms on the right-hand side of (4.4) 
as follows. 
\begin{align}
\begin{aligned}
& \left| \, \int^{\infty}_{-\infty} 
\big( \, \partial_t\phi-\partial_x^2 \phi \, \big) \, 
   \bigl( \, 
    f(\phi +U^r)-f(U^r) 
    \,  \bigr) \, \mathrm{d}x \, \right|\\
& \le \epsilon \left( \, 
\| \, \partial_{t}\phi \, \|_{L^2}^{2}
+ \| \, \partial_x^2\phi \, \|_{L^2}^{2}
\, \right)
+ C_{\epsilon, \beta} \, 
   \left\| \, 
    f(\phi +U^r)-f(U^r) 
    \, \right\|_{L^2}^2,
\end{aligned}
\end{align}
\begin{align}
\begin{aligned}
& \left| \, \int^{\infty}_{-\infty} 
\big( \, \partial_t\phi-\partial_x^2 \phi \, \big) \, 
   F(U^r) \, \mathrm{d}x \, \right|\\
& \le \epsilon \left( \, 
\| \, \partial_{t}\phi \, \|_{L^2}^{2}
+ \| \, \partial_x^2\phi \, \|_{L^2}^{2}
\, \right)
+ C_{\epsilon, \beta} \, 
   \left\| \, 
    F(U^r) 
    \, \right\|_{L^2}^2,
\end{aligned}
\end{align}
for $\epsilon>0$. 
Choosing $\epsilon$ suitably small, substituting (4.5) and (4.6) into (4.4), 
and using Lemma 4.3, we obtain the desired estimate. 

Thus, we complete the proof of Proposition 4.4. 

\medskip

\noindent
{\bf Remark 4.5.}\quad 
Similarly to Lemma 4.2, 
we have the uniform boundedness of $\partial_x\phi$ 
by using Proposition 4.4 that 
\begin{equation*}
\sup_{t\in [0,T], \, x\in \mathbb{R}}\, 
|\, \partial_x\phi(t,x) \, | \le C_{\phi_0}.
\end{equation*}

\medskip

\noindent
We further show the {\it a priori} 
estimate for $\partial_{x}^2\phi$ and $\partial_{x}^3\phi$ as follows. 

\medskip

\noindent
{\bf Proposition 4.6.}\quad {\it
There exists a positive constant $C_{\phi_0}$
such that 
\begin{equation*}
\| \, \partial_{x}^2\phi(t) \, \|_{H^1}^{2} 
  + \int^{t}_{0} 
  \left( \, 
  \| \, \partial_t\partial_x\phi(\tau) \, \|_{H^1}^{2}
+\| \, \partial_{x}^3\phi(\tau) \, \|_{H^1}^{2} 
\, \right)
\, \mathrm{d}\tau
\leq C_{\phi_0} 
\quad \big( \, t \in [\, 0, \, T\, ] \, \big).
\end{equation*}
}

\medskip

{\bf Proof of Proposition 4.6.}
Multiplying the equation in (3.2) by $\partial_x^4\phi-\partial_t\partial_x^2 \phi$ 
and integrating it with respect to $x$, 
we have, after integration by parts, that 
\begin{align}
\begin{aligned}
& \frac{\beta+1}{2} \, \frac{\mathrm{d}}{\mathrm{d}t} \, \| \, \partial_x^2 \phi(t) \, \|_{L^2}^2 
+ \frac{\alpha+\gamma}{2} \, \frac{\mathrm{d}}{\mathrm{d}t} \, \| \, \partial_x^3 \phi(t) \, \|_{L^2}^2 \\
&+ \| \, \partial_{t}\partial_x\phi(t) \, \|_{L^2}^{2} 
+ \alpha \, \| \, \partial_{t}\partial_{x}^2\phi(t) \, \|_{L^2}^{2} 
+ \beta \, \| \, \partial_{x}^3\phi(t) \, \|_{L^2}^{2} 
+ \gamma \, \| \, \partial_{x}^4\phi(t) \, \|_{L^2}^{2} \\
&=-\int^{\infty}_{-\infty} 
\big( \, \partial_x^4\phi-\partial_t\partial_x^2 \phi \, \big) \, 
   \bigl( \, 
    f(\phi +U^r)-f(U^r) 
    \,  \bigr) \, \mathrm{d}x \\
&+ \int^{\infty}_{-\infty} 
\big( \, \partial_x^4\phi-\partial_t\partial_x^2 \phi \, \big) \, F(U^r) \, \mathrm{d}x.
\end{aligned}
\end{align}
The each terms on the right-hand side of (4.7) can be estimated 
quite similarly to (4.5) and (4.6). 
Therefore, we obtain the desired estimate by using Lemma 4.3. 

Thus, we complete the proof of Proposition 4.6. 

\medskip

\noindent
{\bf Remark 4.7.}\quad 
Similarly to Lemma 4.2 and remark 4.5, 
we have the uniform boundedness of $\partial_x^2\phi$ 
by using Proposition 4.6 that 
\begin{equation*}
\sup_{t\in [0,T], \, x\in \mathbb{R}}\, 
|\, \partial_x^2\phi(t,x) \, | \le C_{\phi_0}.
\end{equation*}
\bigskip 

\noindent
\section{Remarks on the uniform estimates}
It is worthwhile to mention the uniform estimates of the solution $\phi$ 
to (3.2).  
By using the Sobolev inequality to 
$\phi$, $\partial_x\phi$, $\partial_x^2\phi$, $\partial_x^3\phi$, 
$\partial_x^3\phi$, $\partial_t\phi$ and $\partial_t\partial_x\phi$, 
and integration by parts, we obtain 
\begin{align}
\begin{aligned}
\sup_{x \in \mathbb{R}}|\,\phi(t,x)\,|
&\le \sqrt{2}\, \| \,\phi(t)\, \|^{\frac{1}{2}}_{L^2} 
     \| \,\partial_x\phi(t)\, \| ^{\frac{1}{2}}_{L^2}\\
&\le \min \left\{ \, 
     \sqrt{2}\, \| \,\phi(t)\, \|^{\frac{3}{4}}_{L^2} 
     \| \,\partial_x^2\phi(t)\, \| ^{\frac{1}{4}}_{L^2}, \, 
     \| \,\phi(t)\, \|_{H^1} \, \right\}\\
&\le \sqrt{2}\, \| \,\phi(t)\, \|^{\frac{3}{4}}_{L^2} 
     \| \,\partial_x\phi(t)\, \| ^{\frac{1}{8}}_{L^2}
     \| \,\partial_x^3\phi(t)\, \| ^{\frac{1}{8}}_{L^2}\\
&\le \sqrt{2}\, \| \,\phi(t)\, \|^{\frac{3}{4}}_{L^2} 
     \| \,\partial_x\phi(t)\, \| ^{\frac{1}{8}}_{L^2}
     \| \,\partial_x^2\phi(t)\, \| ^{\frac{1}{16}}_{L^2}
     \| \,\partial_x^4\phi(t)\, \| ^{\frac{1}{16}}_{L^2}, 
\end{aligned}
\end{align}
\begin{align}
\begin{aligned}
\sup_{x \in \mathbb{R}}|\,\partial_x\phi(t,x)\,|
&\le \sqrt{2}\, \| \,\partial_x\phi(t)\, \|^{\frac{1}{2}}_{L^2} 
     \| \,\partial_x^2\phi(t)\, \| ^{\frac{1}{2}}_{L^2}\\
&\le \min \left\{ \, 
     \sqrt{2}\, \| \,\partial_x\phi(t)\, \|^{\frac{3}{4}}_{L^2}
     \| \,\partial_x^3\phi(t)\, \| ^{\frac{1}{4}}_{L^2}, \, 
     \| \,\partial_x\phi(t)\, \|_{H^1} \, \right\}\\
&\le \sqrt{2}\, \| \,\partial_x\phi(t)\, \|^{\frac{3}{4}}_{L^2} 
     \| \,\partial_x^2\phi(t)\, \| ^{\frac{1}{8}}_{L^2}
     \| \,\partial_x^4\phi(t)\, \| ^{\frac{1}{8}}_{L^2}, 
\end{aligned}
\end{align}
\begin{align}
\begin{aligned}
\sup_{x \in \mathbb{R}}|\,\partial_x^2\phi(t,x)\,|
&\le \sqrt{2}\, \| \,\partial_x^2\phi(t)\, \|^{\frac{1}{2}}_{L^2} 
     \| \,\partial_x^3\phi(t)\, \| ^{\frac{1}{2}}_{L^2}\\
&\le \min \left\{ \, 
     \sqrt{2}\, \| \,\partial_x^2\phi(t)\, \|^{\frac{3}{4}}_{L^2} 
     \| \,\partial_x^4\phi(t)\, \| ^{\frac{1}{4}}_{L^2}, \, 
     \| \,\partial_x^2\phi(t)\, \|_{H^1} \, \right\}, 
\end{aligned}
\end{align}
\begin{equation}
\sup_{x \in \mathbb{R}}|\,\partial_x^3\phi(t,x)\,|
\le \sqrt{2}\, \| \,\partial_x^3\phi(t)\, \|^{\frac{1}{2}}_{L^2} 
     \| \,\partial_x^4\phi(t)\, \| ^{\frac{1}{2}}_{L^2}
\le \| \,\partial_x^3\phi(t)\, \|_{H^1},
\end{equation}
\begin{equation}
\sup_{x \in \mathbb{R}}|\,\partial_t\phi(t,x)\,|
\le \sqrt{2}\, \| \,\partial_t\phi(t)\, \|^{\frac{1}{2}}_{L^2} 
     \| \,\partial_t\partial_x\phi(t)\, \| ^{\frac{1}{2}}_{L^2}
\le \| \,\partial_t\phi(t)\, \|_{H^1},
\end{equation}
\begin{equation}
\sup_{x \in \mathbb{R}}|\,\partial_t\partial_x\phi(t,x)\,|
\le \sqrt{2}\, \| \,\partial_t\partial_x\phi(t)\, \|^{\frac{1}{2}}_{L^2} 
     \| \,\partial_t\partial_x^2\phi(t)\, \| ^{\frac{1}{2}}_{L^2}
\le \| \,\partial_t\partial_x\phi(t)\, \|_{H^1},
\end{equation}
Noting the first term on the {\it a priori} estimates in Theorem 3.3, that is, 
\begin{align}
\begin{aligned}
& \| \, \phi(t) \, \|_{H^3}^{2} 
+ \int^{t}_{0} 
   \big\| \, 
   \big(\, \sqrt{\partial_x U^r} \:  \phi \, \big)(\tau) 
   \, \big\|_{L^2}^{2} \, \mathrm{d}\tau 
  \\
& 
+  \int^{t}_{0} 
\| \, \partial_{x}\phi(\tau) \, \|_{H^3}^{2} \, \mathrm{d}\tau
+ \int^{t}_{0} 
\| \, \partial_{t}\phi(\tau) \, \|_{H^2}^{2} \, \mathrm{d}\tau  
\leq C_{\phi_0},%
\end{aligned}
\end{align}
and using (5.1)-(5.6), we arrive at 
\begin{equation}
\sup_{x \in \mathbb{R}}|\,\phi(t,x)\,|
\le C_{\phi_0}\, 
    \min \left\{ \, 
     \min_{k=1,2,3,4}
     \| \,\partial_x^k\phi(t)\, \|_{L^2}^{\frac{1}{2^k}}, \, 
     \| \,\phi(t)\, \|_{H^1} \, \right\}, 
\end{equation}
\begin{equation}
\sup_{x \in \mathbb{R}}|\,\partial_x\phi(t,x)\,|
\le C_{\phi_0}\, 
    \min \left\{ \, 
     \min_{k=1,2,3}
     \| \,\partial_x^{k+1}\phi(t)\, \|_{L^2}^{\frac{1}{2^k}}, \, 
     \| \,\partial_x\phi(t)\, \|_{L^2}^{\frac{1}{2}}, \, 
     \| \,\partial_x\phi(t)\, \|_{H^1} \, \right\}, 
\end{equation}
\begin{equation}
\sup_{x \in \mathbb{R}}|\,\partial_x^2\phi(t,x)\,|
\le C_{\phi_0}\, 
    \min \left\{ \, 
     \| \,\partial_x^2\phi(t)\, \|_{L^2}^{\frac{1}{2}}, \, 
     \| \,\partial_x^3\phi(t)\, \|_{L^2}^{\frac{1}{2}}, \, 
     \| \,\partial_x^4\phi(t)\, \|_{L^2}^{\frac{1}{4}}, \, 
     \| \,\partial_x^2\phi(t)\, \|_{H^1} \, \right\}, 
\end{equation}
\begin{equation}
\sup_{x \in \mathbb{R}}|\,\partial_x^3\phi(t,x)\,|
\le C_{\phi_0}\, 
    \min \left\{ \, 
     \| \,\partial_x^3\phi(t)\, \|_{L^2}^{\frac{1}{2}}, \, 
     \| \,\partial_x^4\phi(t)\, \|_{L^2}^{\frac{1}{2}}
      \, \right\}, 
\end{equation}
for some $C_{\phi_0}>0$. 
Then, by (5.4), (5.5) and (5.8)-(5.11), 
noting the interpolation, 
we obtain the uniform estimates as follows. 

\medskip

\noindent
{\bf Proposition 5.1.}\quad {\it
There exists a positive constant $C_{\phi_0}$
such that 
\begin{equation*}
\int_0^t
\left( \, 
\sup_{x \in \mathbb{R}}|\,\phi(\tau,x)\,|
\, \right)^k
\, \mathrm{d}\tau
\le C_{\phi_0}
\quad \big( \, 2 \le k \le 32 \, \big), 
\end{equation*}
\begin{equation*}
\int_0^t
\left( \, 
\sup_{x \in \mathbb{R}}|\,\partial_x\phi(\tau,x)\,|
\, \right)^k
\, \mathrm{d}\tau
\le C_{\phi_0}
\quad \big( \, 2 \le k \le 16 \, \big), 
\end{equation*}
\begin{equation*}
\int_0^t
\left( \, 
\sup_{x \in \mathbb{R}}|\,\partial_x^2\phi(\tau,x)\,|
\, \right)^k
\, \mathrm{d}\tau
\le C_{\phi_0}
\quad \big( \, 2 \le k \le 8 \, \big), 
\end{equation*}
\begin{equation*}
\int_0^t
\left( \, 
\sup_{x \in \mathbb{R}}|\,\partial_x^3\phi(\tau,x)\,|
\, \right)^k
\, \mathrm{d}\tau
\le C_{\phi_0}
\quad \big( \, 2 \le k \le 4 \, \big), 
\end{equation*}
\begin{equation*}
\int_0^t
\left( \, 
\sup_{x \in \mathbb{R}}|\,\partial_t\phi(\tau,x)\,|
\, \right)^2
\, \mathrm{d}\tau
\le C_{\phi_0}, 
\end{equation*}
\begin{equation*}
\int_0^t
\left( \, 
\sup_{x \in \mathbb{R}}|\,\partial_t\partial_x\phi(\tau,x)\,|
\, \right)^2
\, \mathrm{d}\tau
\le C_{\phi_0}, 
\end{equation*}
for $t \in [\, 0, \, T\, ]$.
}

\bigskip









\bibliographystyle{model6-num-names}
\bibliography{<your-bib-database>}

\begin{thebibliography}{99}

\bibitem{amik-bona-schonbek}
C.J. Amick, J.L. Bona and M.E. Schonbek, 
{\it Decay of solutions of some nonlinear wave equations}, 
J. Differential Equations {\bf 81} (1989), pp. 1-49.

\bibitem{and-ego-lan-tes}
K. Andreiev, I. Egorova, T.L. Lange and G. Teschl, 
{\it Rarefaction waves of 
the Korteweg-de Vries equation via nonlinear steepest deecent}, 
J. Differential Equations {\bf 10} (2016), pp. 5371-5410.
%

\bibitem{benjamin-bona-mahony}
T.B. Benjamin, J.L. Bona and J.J. Mahony, 
{\it Model equations for long waves in nonlinear dispersive system}, 
Phil. Trans. R. Soc. Lond. Ser. A {\bf 272} (1972), pp. 47-78. 

\bibitem{bona-schonbek}
J.L. Bona and M.E. Schonbek, 
{\it Travelling-wave solutions to the Korteweg-deVries-Burgers equation}, 
Proc. Roy. Soc. Edinburgh {\bf 101A} (1985), pp. 207-226. 

\bibitem{bona-rajopadhye-schonbek}
J.L. Bona, S.V. Rajopadhye and M.E. Schonbek, 
{\it Models for the propagation of bores I. Two dimensional theory}, 
Differential Integral Equations {\bf 7} (1994), pp. 699-734.

\bibitem{chh1}
R.P. Chhabra, 
{\it Bubbles, drops and particles in non-Newtonian Fluids}, 
CRC, Boca Raton, FL (2006).

\bibitem{chh2}
R.P. Chhabra, 
{\it Non-Newtonian Fluids: An Introduction}, 
URL http://www.physics.iitm.ac.in/$\tilde{\: }$compflu/Lect-notes/chhabra.pdf.

\bibitem{chh-ric}
R.P. Chhabra and J.F. Richardson,  
{\it Non-Newtonian flow and applied rheology}, 
2nd edn. Butterworth-Heinemann, Oxford (2008).

\bibitem{de waele} 
A. de Waele, 
{\it Viscometry and plastometry}, 
J. Oil Colour Chem. Assoc. {\bf 6} (1923), pp. 33-69.


\bibitem{dua-fan-kim-xie}
R. Duan, L.-L. Fan, J.-S. Kim and L.-Q. Xie, 
{\it Nonlinear stability of strong rarefaction waves 
for the generalized KdV-Burgers-Kuramoto equation 
with large initial perturbation},
Nonlinear Anal. {\bf 73} (2010), pp. 3254-3267.

\bibitem{dua-zha}
R. Duan and H.-J. Zhao, 
{\it Global stability of strong rarefaction waves 
for the generalized KdV-Burgers equation},
Nonlinear Anal. {\bf 66} (2007), pp. 1100-1117.

\bibitem{ego-gru-tes}
I. Egorova, K. Grunert and G. Teschl, 
{\it On the Cauchy problem for 
the Korteweg-de Vries equation with steplike finite-gap initial data I. 
Schwartz-type perturbations}, 
Nonlinearity {\bf 22} (2009), pp. 1431-1457. 

\bibitem{ego-tes}
I. Egorova and G. Teschl, 
{\it On the Cauchy problem for 
the Korteweg-de Vries equation with steplike finite-gap initial data I\hspace{-.1em}I. 
Perturbations with finite moments}, 
J. Anal. Math. {\bf 115} (2011), pp. 71-101. 


\bibitem{harabetian}
E. Harabetian, 
{\it Rarefaction and large time behavior for parabolic 
equations and monotone schemes}, 
Comm. Math. Phys. {\bf 114} (1988), pp. 527-536.

\bibitem{hashimoto-matsumura} 
I. Hashimoto and A. Matsumura, 
{\it Large time behavior 
of solutions to an initial boundary value problem 
on the half space for scalar viscous conservation law}, 
Methods Appl. Anal. {\bf 14} (2007), pp. 45-59.

\bibitem{hattori-nishihara} 
Y. Hattori and K. Nishihara, 
{\it A note on the stability of rarefaction wave of the 
Burgers equation}, 
Japan J. Indust. Appl. Math. {\bf 8} (1991), pp. 85-96.
%

\bibitem{ilin-kalashnikov-oleinik} 
A.M. Il'in, A.S. Kala{\v{s}}nikov and O.A. Ole{\u\i}nik,
{\it Second-order linear equations of parabolic type}, 
Uspekhi Math. Nauk SSSR {\bf 17} (1962), pp. 3-146 (in Russian);
English translation in Russian Math. Surveys {\bf 17} (1962), pp. 1-143. 

\bibitem{ilin-oleinik} 
A.M. Il'in and O.A. Ole{\u\i}nik, 
{\it Asymptotic behavior of the solutions of the Cauchy problem for 
some quasi-linear equations for large values of the time}, 
Mat. Sb. {\bf 51} (1960), pp. 191-216 (in Russian).

\bibitem{jah-str-mul} 
P. Jahangiri, R. Streblow and D. M\"{u}ller, 
{\it Simulation of Non-Newtonian Fluids using Modelica}, 
Proceedings of the 9th International Modelica Conference September 3-5, 
Munich, Germany, (2012), pp. 57-62. 
%

\bibitem{kanel} 
Y. Kanel', 
{\it On a model system of one-dimensional gas motion}, 
Differencial'nya Uravrenija {\bf 4} (1968), pp. 374-380. 

%
%

\bibitem{kondo-webler1}
C.I. Kondo and C.M. Webler 
{\it The generalized BBM-Burger equations with non-linear 
dissipative term: existence and convergence results}, 
Appl. Anal. {\bf 87} (2008), pp. 977-995.

\bibitem{kondo-webler2}
C.I. Kondo and C.M. Webler 
{\it Higher order for the generalized BBM-Burgers equation: 
existence and convergence results}, 
Appl. Anal. {\bf 88} (2009), pp. 1085-1101.

\bibitem{kondo-webler3}
C.I. Kondo and C.M. Webler 
{\it Higher order for the generalized BBM-Burgers equation: 
existence and convergence results}, 
Acta Appl. Math. {\bf 111} (2010), pp. 45-64.

\bibitem{kondo-webler4}
C.I. Kondo and C.M. Webler 
{\it The generalized BBM-Burgers equations: convergence results 
for conservation law with discontinuous flux function}, 
Appl. Anal. {\bf 95} (2016), pp. 503-523.

\bibitem{lad1} 
O.A. Lady{\v{z}}enskaja, 
{\it New Equations for the Description of the Viscous Incompressible
Fluids and Solvability in the Large of the Boundary Value Problems for Them, in 
``Boundary Value Problems of Mathematical Physics V''}, 
Amer. Math. Soc., 
Providence, Rhode. Island, 1970.

\bibitem{lad-sol-ura} 
O. A. Lady{\v{z}}enskaja, V. A. Solonnikov and N. N. Ural'ceva, 
{\it Linear and quasilinear equations of parabolic type}, 
Transl. Math. Monographs, vol. 23, Amer. Math. Soc., 
Providence, Rhode. Island, 1968.

\bibitem{lax} 
P.D. Lax, 
{\it Hyperbolic systems of conservation laws  I\hspace{-.1em}I}, 
Comm. Pure Appl. Math. {\bf 10} (1957), pp. 537-566.

\bibitem{liep-rosh} 
H.W. Liepmann and A. Roshko, 
{\it Elements of Gas Dynamics}, 
John Wiley \& Sons, Inc., New York, 1957.
%

\bibitem{liu-matsumura-nishihara} 
T.-P. Liu, A. Matsumura and K. Nishihara, 
{\it Behaviors of solutions for the Burgers equation 
with boundary corresponding to rarefaction waves}, 
SIAM J. Math. Anal. {\bf 29} (1998), pp. 293-308.

\bibitem{ma} 
J. M{\'{a}}lek, 
{\it Some frequently used models for non-Newtonian fluids}, 
URL http://www.karlin.mff.cuni.cz/$\tilde{\: }$malek/new/images/Lecture4.pdf.

\bibitem{ma-pr-st} 
J. M{\'{a}}lek, D. Pra{\v{z}}{\'{a}}k and M. Steinhauer,
{\it On the existence and regularity of solutions 
for degenerate power-law fluids}, 
Differential Integral Equations {\bf 19} (2006), pp. 449-462.
%

\bibitem{matsu-nishi1} 
A. Matsumura and K. Nishihara,
{\it Asymptotic toward the rarefaction wave of solutions of a
one-dimensional model system for compressible viscous gas}, Japan
J. Appl. Math. {\bf 3} (1986), pp. 1-13.

\bibitem{matsu-nishi2} 
A. Matsumura and K. Nishihara, 
{\it Asymptotics toward the rarefaction wave 
of the solutions of Burgers' equation with 
nonlinear degenerate viscosity
}
Nonlinear Anal. TMA {\bf 23} (1994), pp. 605-614.

\bibitem{matsu-nishi3} 
A. Matsumura and K. Nishihara, 
{\it Asymptotic stability of traveling waves for scalar 
viscous conservation laws with non-convex nonlinearity}, 
Comm. Math. Phys. {\bf 165} (1994), pp. 83-96.

\bibitem{matsumura-yoshida} 
A. Matsumura and N. Yoshida, 
{\it Asymptotic behavior of solutions to the Cauchy problem 
for the scalar viscous conservation law 
with partially linearly degenerate flux}, 
SIAM J. Math. Anal. {\bf 44} (2012), pp. 2526-2544. 

\bibitem{matsumura-yoshida'} 
A. Matsumura and N. Yoshida, 
{\it Global asymptotics toward 
the rarefaction waves for solutions to the Cauchy problem of
the scalar conservation law with nonlinear viscosity}, 
Osaka J. Math. {\bf 57} (2020), pp. 187-205.

\bibitem{mei1}
M. Mei, 
{\it Large-time behavior of solution for 
generalized Benjamin-Bona-Mahony-Burgers equations},
Nonlinear Anal. {\bf 33} (1998), pp. 699-714.

\bibitem{mei2}
M. Mei, 
{\it $L^q$-decay rates of solutions for 
generalized Benjamin-Bona-Mahony-Burgers equations},
J. Differential Equations {\bf 158} (1999), pp. 314-340.

\bibitem{mei-schmeiser}
M. Mei and C. Schmeiser, 
{\it Asymptotic profiles of solutions for the BBM-Burgers equation},
Funkcialaj Ekvacioj {\bf 44} (2001), pp. 151-170.

\bibitem{naumkin}
P.I. Naumkin, 
{\it Large-time asymptotic of a step for the Benjamin-Bona-Mahony-Burgers equation}, 
Proc. Roy. Soc. Edinburgh {\bf 126A} (1996), pp. 1-18. 

\bibitem{nishi-raj}
K. Nishihara and S.V. Rajopadhye, 
{\it Asymptotic behavior of solutions to the Korteweg-deVries-Burgers equation}, 
Differential Integral Equations {\bf 11} (1998), pp. 85-93.
%
%

\bibitem{osh-ral} 
S. Osher and J. Ralston
{\it $L^1$ stability of traveling waves with applications to convective porous media flow}, 
Comm. Pure Appl. Math. {\bf 35} (1982), pp. 737-751.

\bibitem{ost} 
W. Ostwald, 
{\it \"{U}ber die Geschwindigkeitsfunktion der Viskositat disperser Systeme}, 
I. Colloid Polym. Sci. {\bf 36} (1925), pp. 99-117 (in German).
%

\bibitem{peregrine}
D.H. Peregrine, 
{\it Calculations of the development of an undular bore}, 
J. Fluid Mech. {\bf 25} (1966), pp. 321-330. 

\bibitem{rajopadhye}
S.V. Rajopadhye, 
{\it Decay rates for the solutions of model equations for bore propagation}, 
Proc. Roy. Soc. Edinburgh {\bf 125A} (1995), pp. 371-398. 

\bibitem{ras-nik-kho}
J. Rashindinia, O. Nikan and L. Khoddam, 
{\it Numerically stable scheme to approximate the nonlinear
KdV-Benjamin-Bona-Mahony-Burger's equation},
4th International Conference on Combinatorics, Cryptography, Computer Science and Computing (2019), pp. 1-10.

\bibitem{rua-gao-che}
L.-Z. Ruan, W.-L. Gao and J. Chen, 
{\it Asymptotic stability of the rarefaction wave 
for the generalized KdV-Burgers-Kuramoto equation},
Nonlinear Anal. {\bf 68} (2008), pp. 402-411.


\bibitem{soc} 
T. Sochi, 
{\it Pore-Scale Modeling of Non-Newtonian Flow in Porous Media}, 
PhD thesis, Imperial College London, 2007. 
%
%

\bibitem{wang} 
Y. Wang, 
{\it On time periodic solutions to the generalized 
BBM-Burgers equation with time-dependent periodic external force},
Math. Model. Anal. {\bf 25} (2020), pp. 184-197.

\bibitem{wan-zhu} 
Z.-A. Wang and C.-J. Zhu, 
{\it Stability of the rarefaction wave 
for the generalized KdV-Burgers equation},
Acta Math. Sci. {\bf 22B}(3) (2002), pp. 319-328.

\bibitem{xu-li} 
H. Xu and B. Li, 
{\it Global existence and bounded estimate of solutions of the BBM-Burgers equation},
Wuhan Univ. J. Nat. Sci. {\bf 21} (2016), pp. 428-432.

\bibitem{yin-zhao-kim}
H. Yin, H. Zhao and J. Kim, 
{\it Convergence rates of solutions toward boundary layer solutions  
for generalized Benjamin-Bona-Mahony-Burgers equations in the half-space}, 
J. Differential Equations {\bf 245} (2008), pp. 3144-3216.

\bibitem{yoshida1} 
N. Yoshida, 
{\it Decay properties of solutions toward a multiwave pattern 
for the scalar viscous conservation law 
with partially linearly degenerate flux}, 
Nonlinear Anal. {\bf 96} (2014), pp. 189-210.

\bibitem{yoshida2} 
N. Yoshida, 
{\it Decay properties of solutions 
to the Cauchy problem for the scalar conservation law 
with nonlinearly degenerate viscosity}, 
Nonlinear Anal. {\bf 128} (2015), pp. 48-76.

\bibitem{yoshida3'} 
N. Yoshida, 
{\it Large time behavior of solutions toward a multiwave pattern for the Cauchy problem of the 
scalar conservation law with degenerate flux and viscosity}, 
in: Mathematical Analysis in Fluid and Gas Dynamics, 
S\={u}rikaisekikenky\={u}sho K\={o}ky\={u}roku {\bf 1947} (2015), pp. 205-222.

\bibitem{yoshida3} 
N. Yoshida, 
{\it Asymptotic behavior of solutions toward a multiwave pattern 
for the scalar conservation law 
with the Ostwald-de Waele-type viscosity}, 
SIAM J. Math. Anal. {\bf 49} (2017), pp. 2009-2036.

\bibitem{yoshida4} 
N. Yoshida, 
{\it Decay properties of solutions toward a multiwave pattern 
to the Cauchy problem for the scalar conservation law 
with degenerate flux and viscosity}, 
J. Differential Equations {\bf 263} (2017), pp. 7513-7558.

\bibitem{yoshida5} 
N. Yoshida, 
{\it Asymptotic behavior of solutions toward the viscous shock waves 
to the Cauchy problem 
for the scalar conservation law with nonlinear flux and viscosity}, 
SIAM J. Math. Anal. {\bf 50} (2018), pp. 891-932.

\bibitem{yoshida6} 
N. Yoshida, 
{\it Asymptotic behavior of solutions toward 
the rarefaction waves to the Cauchy problem for 
the scalar diffusive dispersive conservation laws}, 
Nonlinear Anal. {\bf 189} (2019), pp. 1-19.

\bibitem{yoshida7} 
N. Yoshida, 
{\it Global structure of solutions toward 
the rarefaction waves for the Cauchy problem 
of the scalar conservation law with nonlinear viscosity}, 
J. Differential Equations {\bf 269} (2020), pp. 10350-10394.

\bibitem{yoshida8} 
N. Yoshida, 
{\it Asymptotic behavior of solutions toward a multiwave pattern 
to the Cauchy problem for the dissipative wave equation 
with partially linearly degenerate flux}, 
Funkcialaj Ekvacioj {\bf 64} (2021), pp. 49-73.

\bibitem{yoshida9} 
N. Yoshida, 
{\it Asymptotic behavior of solutions toward the constant state to the Cauchy problem for the 
non-viscous diffusive dispersive conservation law}, 
Preprint, arXiv:2107.07874. 

\bibitem{yoshida10} 
N. Yoshida, 
{\it Global asymptotic stability of the rarefaction waves to the Cauchy problem for the scalar 
non-viscous diffusive dispersive conservation laws}, 
Preprint.

\bibitem{zhao-xuan} 
H. Zhao and B. Xuan, 
{\it Existence and convergence of solutions 
for the generalized BBM-Burgers equations with dissipative term}, 
Nonlinear Anal. {\bf 28} (1997), pp. 1835-1849.

%

\end{thebibliography}







\end{document}